\documentclass{article}
\usepackage{amsmath,amsthm,amssymb}
\usepackage[a4paper]{geometry}

\newtheorem{theorem}{Theorem}[section]

\newtheorem{proposition}[theorem]{Proposition}

\theoremstyle{definition}
\newtheorem{definition}[theorem]{Definition}
\newtheorem*{remark}{Remark}
\newtheorem{example}[theorem]{Example}

\newcommand{\R}{\mathbb{R}}
\newcommand{\C}{\mathbb{C}}
\newcommand{\Z}{\mathbb{Z}}
\newcommand{\F}{\mathcal{F}}
\newcommand{\E}{\varepsilon}

\title{A note on Mitsumatsu's construction of\\
a leafwise symplectic foliation}\author{Atsuhide Mori}
\date{}
\begin{document}

\maketitle

\begin{abstract}
Mitsumatsu \cite{Mi} constructed 
a leafwise symplectic structure of the Lawson foliation of $S^5$. 
Following his construction, we improve a 
previous result of the author \cite{Mo}
on convergence of contact structure to codimension one foliation,
and give a sufficient condition 
for convergence to foliation with leafwise symplectic structure. 
As an application, we show that the product $S^4\times S^1$ admits 
infinitely many codimension one leafwise symplectic foliations.
\end{abstract}

\maketitle

\section{Introduction}
After efforts of Verjovsky et al., Mitsumatsu \cite{Mi} 
succeeded in constructing a leafwise symplectic structure 
of the Lawson foliation of $S^5$. 
His construction is based on an open-book decomposition of $S^5$ 
whose binding is a $T^2$-bundle over $S^1$. 
It also works on infinitely many foliations of $S^5$
which appear in the author's  
result \cite{Mo} (and its completion by Kasuya \cite{Kas}) 
on convergence of a contact structure to a foliation with dead-end components. 
It is natural to ask if such convergence 
has anything to do with leafwise symplecticity in general.  

The pair $(\alpha,\tau)$ 
of a $1$-form $\alpha$ and a $2$-form $\tau$ on an oriented 
$(2n+1)$-manifold is called a twisted contact structure 
if it satisfies $\alpha\wedge(d\alpha+\tau)^n>0$ 
(Nunes da Costa and Petalidou \cite{NP}).
Suppose that for any $\E\in (0,1]$, 
the pair $(\alpha, \E\tau)$ is twisted contact. 
Then we call $\alpha$ an $\E\tau$-confoliation form 
since it satisfies $\alpha\wedge(d\alpha)^n\ge 0$. 
Eliashberg and Thurston \cite{ET} introduced 
the notion of confoliation 
and developed it in three dimension.
Higher dimensional confoliation was developed by Altschuler and Wu \cite{AW}
under an additional ``conductivity'' condition which ironically foliations never satisfy. 
In this article, we connect a contact structure to  
a leafwise symplectic foliation   
through $\E\tau$-confoliations. 
The key observation is the following proposition. 

\begin{proposition}[Proposition~\ref{prop}]
Suppose that we have a contact form $\eta$, 
a non-singular closed $1$-form $\nu$, and 
a $2$-form $\tau$ on a closed $(2n-1)$-manifold 
such that $\eta_t=(1-t)\nu+t\eta$ is 
$\E\tau$-confoliation for each $t\in [0,1]$. 
Then if $n>1$, we have $\nu\wedge (d\eta)^{n-1}=0$ and $\nu\wedge \tau^{n-1}>0$. 
Further the $1$-form $\eta_t$ is contact for $t>0$. 
\end{proposition}

We abstract the essence of the construction of Mitsumatsu
to show the following extension theorem. 

\begin{theorem}[Theorem~\ref{thm}]
Let $\alpha$ be a contact form on a closed oriented $(2n+1)$-manifold $M^{2n+1}$ 
and $\Sigma_0$ a page of an open-book decomposition of $M^{2n+1}$ which 
supports the contact structure $\ker\alpha$. 
Suppose that the binding $N^{2n-1}=\partial\Sigma$ admits 
a family $\{\eta_t\}_{t\in[0,1]}$ of confoliation forms 
described in the above proposition such that 
$\eta_1=\alpha|_{N^{2n-1}}$. Precisely, we assume that 
$\eta_t=(1-t)\nu+t\eta_1$ is $\E\tau_N$-confoliation 
on the binding $N^{2n-1}$ for each $t\in[0,1]$, 
where $\nu$ is a non-singular closed $1$-form and $\tau_N$ is 
a $2$-form. Then it extends to a family $\{\alpha_t\}$ 
of $\E\tau$-confoliation forms on $M^{2n+1}$ with respect to 
a $2$-form $\tau$ such that 
$\ker\alpha_0$ is the tangent space of a regular foliation $\F$,  
$\{\ker \alpha_t\}_{t \in (0,1]}$ is a family of contact structures 
supported by the open-book decomposition, and  
$\ker\alpha_1$ is the original contact structure $\ker\alpha$. 
Further if $\tau_N$ is closed and extends to a closed $2$-form on the page $\Sigma_0$, 
we can take $\tau$ so that $(d\alpha_0+\E\tau)|_{T\F}=\E\tau|_{T\F}$ is a leafwise symplectic form. 
\end{theorem}

This provides a new example of leafwise symplectic foliation of $S^4\times S^1$, 
which can be deformed into the contact structure on the strongly pseudo convex 
boundary of $B^5\times S^1$ (\S 4). 

\section{Definitions and statement of result}
Let $\alpha$ be a $1$-form and $\tau$ a $2$-form on 
a closed oriented $(2n+1)$-manifold $M^{2n+1}$. 

\begin{definition} 
We call the pair $(\alpha,\tau)$ a {\em twisted contact structure} if 
it satisfies $\alpha\wedge(d\alpha+\tau)^n>0$. 
Then $\tau$ is called the {\em twisting}. 
We call $\alpha$ an {\em $\E\tau$-confoliation form} and $\tau$ its {\em reference twisting} 
if the pair $(\alpha, \E\tau)$ is a twisted contact structure for any $\E\in(0,1]$. 
\label{tau}
\end{definition}
\begin{remark} 
If $\alpha$ is an $\E\tau$-confoliation form with reference twisting $\tau$, 
then clearly $e^f\alpha$ is an $\E\tau'$-confoliation form with reference twisting $\tau'=e^f\tau$ for any $f\in C^\infty(M^{2n+1})$. 
Note that the difference between $d(e^f\alpha)$ and $e^fd\alpha$ 
vanishes along $\ker\alpha$, i.e., 
$\alpha\wedge d(e^f\alpha)=e^f\alpha\wedge d\alpha$ holds. 
Further, if we define $(\E\tau, h)$-confoliation form by replacing the interval 
$(0,1]$ for $\E$ 
with the interval $(0,h]$ of arbitrary length $h>0$, then we can see that for 
any (compactly supported) functions $f$ and $g$, 
there exists a length $h>0$ such that $e^f\alpha$ is an $(\E\tau'', h)$-confoliation form 
with reference twisting $\tau''=e^g\tau$. 
In other words, we fix the scale of the Pfaff form $\alpha$ of the hyperplane field 
and that of the $2$-form $\tau$ in order to 
set the length $h$ of the deformation to the unit.  
\end{remark}
\begin{example}
A hyperplane field $\ker\alpha$ on $M^{2n+1}$ 
is a contact structure if and only if $(\alpha,0)$ is a twisted contact structure. 
On the other hand, in the case where 
$\ker\alpha$ defines a codimension one foliation (i.e. $\alpha\wedge d\alpha=0$),
the restriction $\tau|_{\ker\alpha}$ defines a leafwise almost symplectic structure 
if and only if $(\alpha,\tau)$ is a twisted contact structure. 
In $1$-dimensional case (i.e. $n=0$, $M^1=S^1$), these examples are merged together. 
Namely, $\tau=0$ is a leafwise symplectic structure of 
the point foliation defined by a contact $1$-form $\alpha>0$. 
\end{example} 
The name of $\E\tau$-confoliation comes from 
the following relevant notions on a hyperplane field $\ker\alpha$. 
\begin{definition}
\begin{itemize}
\item[(i)] Suppose that an almost complex structure $J$ on 
$\ker\alpha$ satisfies $d\alpha(v, J v)\ge0$ 
for any $v\in\ker\alpha$. 
Then the pair $(\ker\alpha,J)$ is called a {\em pseudo-convex 
almost CR structure} on $M^{2n+1}$. 
Eliashberg and Thurston \cite{ET} studied the topology of  
a hyperplane field $\ker\alpha$ which admits 
an almost complex structure $J$  
such that $(\ker\alpha,J)$ is a pseudo convex almost CR structure. 
We call such a hyperplane field an {\em ET-confoliation}.
We should notice that they actually studied it 
only in $3$-dimensional case.

\item[(ii)] Altschuler and Wu \cite{AW} considered a hyperplane field 
$\ker\alpha$ ($\alpha\neq 0$) which satisfies $\alpha\wedge(d\alpha)^n\ge 0$. 
We call it an {\em AW-confoliation}.  
They studied it in dimension $\ge 3$ under the following condition called 
{\em conductivity}. The $3$-dimensional conductivity is the accessibility of 
each degenerate point to a non-degenerate (i.e. contact) point 
by an arc tangent to the AW-confoliation. 
The high dimensional conductivity is basically the local splitting property into 
a $(2n-2)$-dimensional symplectic structure and a $3$-dimensional 
conductive AW-confoliation which are defined 
by restricting respectively $d \alpha$ and $\ker\alpha$ 
(see \cite{AW} for the precise definition). 
This implies at least the next-highest non-integrability, i.e., $\alpha\wedge (d\alpha)^{n-1}\neq0$. 

\item[(iii)] A relevant notion due to Massot, Niederkr\"uger, and Wendl \cite{Massot}:  
Let $\alpha$ be a contact form on the boundary $M^{2n+1}=\partial W^{2n+2}$ 
of a compact symplectic manifold $(W^{2n+2},\omega)$, Suppose that 
$\alpha\wedge(td\alpha+\omega|_{M^{2n+1}})^n>0$ holds for any $t\ge 0$. 
Then we say that $(W^{2n+2},\omega)$
is a {\em weak symplectic filling} of $(M, \ker\alpha)$. 
\end{itemize}\noindent
We notice that an AW-confoliation as well as an ET-confoliation 
is simply called a {\em confoliation} in literatures.  
\end{definition}
\begin{remark} The above notions 
relate to $\E\tau$-confoliation as follows.
\begin{itemize}
\item[(i)]
Let $\alpha$ be an $\E\tau$-confoliation form. 
Suppose that an 
almost complex structure $J$ on $\ker\alpha$ 
satisfies the inequality $(d\alpha+\E\tau)(v, J v)>0$ 
for any $v(\in\ker\alpha)\neq 0$ and any $\E\in(0,1]$. 
Then the $\E\tau$-confoliation form $\alpha$ defines 
an ET-confoliation 
since the inequality becomes $d\alpha(v, J v)\geq 0$ as $\E \to 0$.
\item[(ii)] 
We see that 
any $\E\tau$-confoliation form defines an AW-confoliation 
by getting $\E$ in Definition~\ref{tau} to go to $0$. 
The converse does not hold 
in the case where $n>1$ 
since any foliation is an AW-confoliation. 
\item[(iii)] Let $\ker\alpha$ be a hyperplane field 
on the boundary $M^{2n+1}=\partial W^{2n+2}$ 
of a compact symplectic manifold $(W^{2n+2}, \omega)$, 
Put $\tau=\omega|_{M^{2n+1}}-d\alpha$. 
Then we see that the condition 
$\alpha\wedge(td\alpha+\omega|_{M^{2n+1}})^n>0$ ($\forall t\ge 0$) is equivalent to 
$\alpha\wedge(d\alpha+\E\tau)^n>0$ ($\forall \E\in(0,1]$) 
which appears in Definition~\ref{tau}. 
If moreover $\alpha$ is contact, $(W^{2n+2},\omega)$ is 
a weak symplectic filling of $(M^{2n+1},\ker\alpha)$. 
\end{itemize}
It was shown in \cite{Massot} that a compact symplectic manifold 
$(W^{2n+2},\omega)$ is a weak symplectic filling of a contact manifold 
$(M^{2n+1}=\partial W^{2n+1},\alpha)$ 
if and only if it admits an almost complex structure $J$ 
such that $\omega(u,Ju)>0$ ($\forall u(\in TW^{2n+2})\neq 0$) holds, 
$\ker\alpha$ is $J$-invariant, and $d\alpha(v,Jv)>0$ 
($\forall v(\in \ker\alpha)\neq 0$) holds. 
Thus the boundary condition of a weak symplectic filling 
provides a weak symplectic version of strong pseudo-convexity. 
From this point of view, our notion of $\E\tau$-confoliation 
is a natural weak symplectic version of pseudo-convexity. 
\end{remark}
In the case where $n>0$, 
it was shown in \cite{AW} that a geometric heat flow 
on the space of $1$-forms drifts 
a conductive AW-confoliation into a contact structure. 
However, there is no foliation satisfying the conductivity condition. 
Our advantage is that a defining $1$-form of 
a foliation $\F$ with a leafwise almost symplectic structure $\tau|_{T\F}$
is an $\E\tau$-confoliation form. 
Indeed the following proposition and example provide a path 
of $\E\tau$-confoliation forms which connects a contact structure 
to such a foliation through contact structures.  

\begin{proposition}
Suppose that a closed $(2n+1)$-manifold $M^{2n+1}$ admits 
a contact form $\alpha$, a non-singular closed $1$-form $\nu$, and 
a $2$-form $\tau$ 
such that the family $\{(1-t)\nu+t\alpha\}_{t\in[0,1]}$ consists of 
$\E\tau$-confoliation forms. 
Then if $n>0$ we have $\nu\wedge (d\alpha)^n=0$ and $\nu\wedge \tau^n>0$. 
Further the $1$-form $(1-t)\nu+t\alpha$ 
is contact for $t>0$. 
\label{prop}
\end{proposition}

\begin{proof} 
The $(2n+1)$-form $\nu\wedge (d\alpha)^n$ 
is null-cohomologous since it is the exterior derivative 
of $\alpha\wedge \nu\wedge (d\alpha)^{n-1}$ ($n>0$). 
On the other hand, it is non-negative 
since the assumption 
\[
\{(1-t)\nu+t\alpha\}\wedge (td\alpha+\E\tau)^n>0
\quad (0\le t\le 1,\quad 0<\E\le 1)
\] 
implies $\{(1-t)\nu+t\alpha\}\wedge (d\alpha)^n\ge 0$ (as $\E\to 0$). 
Thus $\nu\wedge (d\alpha)^n$ vanishes everywhere. 
The assumption also implies $\nu\wedge \tau^n>0$ (when $t=0$). 
Thus $\{(1-t)\nu+t\alpha\}\wedge(td\alpha)^n=
t^{n+1}\alpha\wedge(d\alpha)^n>0$ holds for $t>0$. 
\end{proof}
Then the manifold $M^{2n+1}$ is regularly fibered by hypersurfaces 
so that the fibration approximates the Riemannian foliation $\F$ defined 
by the closed $1$-form $\tau$ at $t=0$ (see Sacksteder \cite{Sa} and Tischler \cite{Ti}) . 
The leafwise almost symplectic structure $\tau|_{T\F}$ is 
leafwise symplectic in the case where   
$\nu\wedge d\tau=0$ holds. 
If moreover $\tau$ is closed, 
the pair $(\nu,\tau)$ defines 
a cosymplectic structure, i.e., satisfies $d\nu=0$, $d\tau=0$, 
and $\nu\wedge \tau^n>0$.
\begin{example}
Let $\beta$ be a contact form 
on a closed oriented $(2n-1)$-manifold $A^{2n-1}$ 
and $\gamma$ a contact form on 
the orientation-reversal $-A^{2n-1}$. 
Suppose that, on the product $\R\times A^{2n-1}$ ($s\in\R$), 
the $2$-form $\omega=d(e^s\beta+e^{-s}\gamma)$ is symplectic. 
Then we say that $\beta$ and $\gamma$ form a Liouville pair on $A^{2n-1}$. 
Such $A^{2n-1}$ exists for any $n>0$. Indeed 
Massot et. al. \cite{Massot} constructed a Liouville pair 
on cusp cross sections of Hilbert modular varieties.  
Then on the closed product manifold $M^{2n+1}=A^{2n-1}\times T^2$ 
($T^2=(\R/2\pi\Z)^2\ni(\theta,\varphi)$), for each $m\in\Z_{>0}$, we put 
\begin{align*}
\Omega&= \left(\frac{1-\cos(m\varphi)}{2}\beta
-\frac{1+\cos(m\varphi)}{2}\gamma \right)
\wedge \left(\frac{1-\cos(m\varphi)}{2}d\beta
+\frac{1+\cos(m\varphi)}{2}d\gamma \right)^{n-1}
\wedge d\theta\wedge d\varphi,\\
\alpha&= \left(\frac{1-\cos(m\varphi)}{2}\beta
+\frac{1+\cos(m\varphi)}{2}\gamma \right)
+\sin(m\varphi) d\theta,\\
\nu&=d\varphi,
\end{align*}
and
\[
\tau
=
\left(\frac{1-\cos(m\varphi)}{2}\beta
-\frac{1+\cos(m\varphi)}{2}\gamma \right)\wedge d\theta
+ \left(\frac{1-\cos(m\varphi)}{2}d\beta
+\frac{1+\cos(m\varphi)}{2}d\gamma \right)
~(\neq~0).
\]
Since the partition $\displaystyle \frac{1\pm\cos(m\varphi)}{2}\left(=\cos^2\frac{m\varphi}{2},\quad \sin^2\frac{m\varphi}{2}\right)$ 
of unity in the expression of $\Omega$ presents the partition $\displaystyle\frac{e^{\pm s}}{e^s+e^{-s}}$ in the following volume form 
on $\R\times A^{2n-1}$, we can see that $\Omega$ is a positive volume form on $M^{2n+1}$. 
\[
\left(\frac{1}{e^s+e^{-s}}\omega\right)^n=ds\wedge 
\left(\frac{e^s}{e^s+e^{-s}}\beta
-\frac{e^{-s}}{e^s+e^{-s}}\gamma \right)
\wedge \left(\frac{e^s}{e^s+e^{-s}}d\beta
+\frac{e^{-s}}{e^s+e^{-s}}d\gamma \right)^{n-1}
\] 
Although the following calculation looks complicated, indeed it is  
almost the same as that in the case where $n=1$. 
From it we see that the assumptions of Proposition~\ref{prop} are satisfied.
\begin{align*}
d\alpha
&=
\frac{m\sin(m\varphi)}{2}d\varphi\wedge(\beta-\gamma)
-m\cos(m\varphi) d\theta\wedge d\varphi
+ \left(\frac{1-\cos(m\varphi)}{2}d\beta
+\frac{1+\cos(m\varphi)}{2}d\gamma \right),
\\
\alpha\wedge (d\alpha)^n
&=mn \left\{-\cos(m\varphi) \left(
\frac{1-\cos(m\varphi)}{2}\beta+\frac{1+\cos(m\varphi)}{2}\gamma
 \right)+\frac{\sin^2(m\varphi)}{2}(\beta-\gamma)
 \right\}\wedge d\theta\wedge d\varphi
\\
&\qquad
\wedge \left(\frac{1-\cos(m\varphi)}{2}d\beta
+\frac{1+\cos(m\varphi)}{2}d\gamma \right)^{n-1}
\\&=
 m n \Omega >0, 
\\
\nu\wedge (t d\alpha+\E\tau )^n
&= n \E d\varphi\wedge \left(\frac{1-\cos(m\varphi)}{2}\beta
-\frac{1+\cos(m\varphi)}{2}\gamma \right)\wedge d\theta
\\
&\qquad
\wedge \left\{(t+\E) \left(\frac{1-\cos(m\varphi)}{2}d\beta
+\frac{1+\cos(m\varphi)}{2}d\gamma \right) \right\}^{n-1}
\\
&= n \E  (t+\E  )^{n-1} \Omega >0\quad(\Rightarrow \nu\wedge (d\alpha)^n=0), 
\\
\alpha\wedge\tau^n
&=0,
\end{align*}
and
\begin{align*}
\alpha\wedge(t d\alpha+\E \tau)^n
&= nt\alpha\wedge d\alpha
\wedge \left\{(t+\E) \left(\frac{1-\cos(m\varphi)}{2}d\beta
+\frac{1+\cos(m\varphi)}{2}d\gamma \right) \right\}^{n-1} \\
&= m n t  (t+ \E  )^{n-1} \Omega \ge 0 \quad (t\ge 0).
\end{align*}
This model is an improvement of Lutz tube in \cite{Mo}
or Giroux torsion domain in \cite{Massot}. 
In the case where $n=1$, we have $M^3=T^3\ni(s,\theta,\varphi)$ with
$\beta=-\gamma=ds$, $\alpha=-\cos(m\varphi)ds+\sin(m\varphi)d\theta$, $\nu=d\varphi$ and 
$\tau=ds\wedge d\theta$. Then the contact structure $\ker\alpha$ is 
called a propeller with respect to the tori $\{\varphi=\mathrm{const}\}$. 
We can imagine that the vertical axis $\partial_\varphi$ of the propeller 
has come up from a horizontal direction via   
$t\partial\varphi+(1-t)\{\cos(m\varphi)\partial_s-\sin(m\varphi)\partial_\theta\}$ like a tippe top. 
This presents how the foliation by the horizontal tori deforms 
to the vertical propeller. 
\end{example}

To obtain topologically more interesting confoliations, 
we move from product manifolds 
to the following setting of open-book decomposition. 
Let $M^{2n+1}$ be a closed oriented $(2n+1)$-manifold, 
$N^{2n-1}\subset M^{2n+1}$ 
a codimension two closed submanifold with trivial neighborhood
$N^{2n-1}\times D^2$, and  
$(\sqrt{\rho}, \theta)$ the polar coordinates of $D^2$. 
Suppose that $\theta$ extends to a fiber bundle projection 
$\theta: M^{2n+1}\setminus N^{2n-1} \to S^1$. 
Let $\Sigma_\theta$ denote the fibers
$\{\theta=\mathrm{const}\}$ which we call the 
pages of an open-book decomposition. 
Then we call $N^{2n-1}$ their binding.
Suppose that there exists a contact form $\alpha$ on $M^{2n+1}$ satisfying
$\alpha\wedge (d\alpha)^{n-1}|_{N^{2n-1}}>0$ and 
$(d\alpha)^n|_{\Sigma_\theta}>0\,\,(\forall \theta)$. 
Let $\eta$ denote the contact form $\alpha|_{N^{2n-1}}$ on the binding $N^{2n-1}$. 
Then modifying $\alpha$ under the condition 
\[
d\theta\wedge (d\alpha |_{M^{2n+1}\setminus N^{2n-1}})^n>0 
\]
if necessary, we can write   
\[
\alpha|_{N^{2n-1}\times D^2}=f(\rho)\eta+g(\rho)d\theta, 
\]
where $\{(f(\rho),g(\rho))\mid 0\le\rho\le 1\}$ is 
a parametrized curve on $\R^2$ such that
\begin{itemize} 
\item[(i)] the initial velocity is $(0,1)$ near the starting point $(f(0),g(0))=(1,0)$, 
\item[(ii)] the final velocity is $(-1,0)$ near the ending point $(f(1),g(1))=(1/2,1)$, 
\item[(iii)] the angular momentum  $L(\rho)=f(\rho)g'(\rho)-f'(\rho)g(\rho)$ is identically $1$, and
\item[(iv)] $f'(\rho)<0$ and $g'(\rho)>0$ hold after $f$ departs from $1$ 
and before $g$ reaches to $1$.
\end{itemize}
Note that both of the area of the triangle 
with vertices $(1,0)$, $(1/2,1)$, and $(0,0)$, 
and that of the area of the fan-shaped region for the above 
curve are equal to $1/2$. 
Giroux \cite{Giroux} proved that 
any closed contact manifold 
of any odd dimension
can be described as above.
Then we say that the contact structure 
is \textit{supported} by the open-book decomposition. 
In the case where a contact form is already specified, 
we can \textit{adapt} it to the above conditions 
by multiplying to it a suitable positive function. 
Now we can state our main result as follows.

\begin{theorem}
Let $M^{2n+1}$, $N^{2n-1}$, $\theta$, $\Sigma_\theta$, 
$\alpha$, and $\eta$ be as above. Namely, we start with the standard description of a contact manifold 
by means of a supporting open-book decomposition and assume that 
the contact form $\alpha$ is adapted. 
Suppose that there exist a non-singular closed $1$-form $\nu$
and a $2$-form $\tau_N$ on $N^{2n-1}$  
such that the family $\{(1-t)\nu+t\eta\}_{t\in[0,1]}$ consists of 
$\E\tau_N$-confoliation forms on $N^{2n-1}$, namely, it satisfies the condition 
described in Proposition~\ref{prop}. 
Then it extends to a family $\{\alpha_t\}_{t\in[0,1]}$ of 
$\E\tau$-confoliation forms on $M^{2n+1}$ 
with respect to a ($t$-independent) $2$-form $\tau$ such that
\begin{itemize} 
\item[(i)] initially at $t=0$, the $\E\tau$-confoliation $\alpha_0$ defines a codimension one foliation 
$\F$, and 
\item[(ii)] after that $\alpha_t$ ($0<t\le 1$) is a contact form which finally becomes 
$\alpha_1=\alpha$ at $t=1$. (Then Gray's stability theorem \cite{Gray}
implies that the contact structure $\ker\alpha_t$ is isotopic to 
$\ker\alpha$ for any $t\in(0,1]$.)   
\end{itemize}
In the case where $\tau_N$ itself is closed and further it extends to a closed $2$-form 
on the page $\Sigma_0$, we can take $\tau$ so that the restriction $\tau|_{T\F}$ is a 
leafwise symplectic form, while $\tau$ itself can not be closed.
\label{thm}
\end{theorem}

\begin{remark}  
Mitsumatsu \cite{Mi} showed that the $2$-form of 
the natural cosymplectic structure by tori on the link of the singularity $(0,0,0)$ 
of the complex surface $z_1^3+z_2^3+z_3^3(+z_1z_2z_3)=0$ 
can be extended to a 
closed $2$-form $\omega$ on the Milnor fiber $\Sigma_0$ 
by means of elementary homology calculation. (In \S 4.3 we give $\omega$ explicitly.) 
This is one of the key steps in his construction of a leafwise symplectic structure 
of the Lawson foliation. 
\end{remark}

\section{Proof of Theorem~\ref{thm}}
Let $\omega$ be any $2$-form on the manifold $M^{2n+1}$ whose restriction to 
the tube $N^{2n-1}\times D^2$ is the pull-back of $\tau_N$ 
under the natural projection to the first factor. 
Then the assumption of the theorem says that 
\begin{equation}
((1-t)\nu+t\eta)\wedge(t d\eta+\E \omega)^{n-1}
\wedge d\rho \wedge d\theta>0\quad (\forall t\in[0,1], \forall \E\in(0,1])
\label{ineq}
\end{equation}
holds on $N^{2n-1}\times D^2$. 
Take a bump function $b(\rho)$ with 
$b(\rho)=1$ for $\rho\le 1/4$, 
$b'(\rho)<0$ for $1/4<\rho<1/2$, 
and $b(\rho)=0$ for $\rho\ge 1/2$. 
We see that the following functions of $\rho\in [0,1]$ are all non-negative 
for any $t\in[0,1]$. 
\begin{align*}
\overline{\,b\,}(\rho) &= \,b\,(1-\rho)~\ge0
&&(\overline{\,b\,}(\rho)=0 ~\Leftrightarrow~ \rho\le1/2),\\
f_t(\rho) &= \{(1-t)\,b\,(\rho)+t\}f(\rho)~\ge0
&&(f_t(\rho)=0 ~\Leftrightarrow~ t=0~\textrm{and}~ \rho\ge 1/2), \\
g_t(\rho) &= \{(1-t)\overline{\,b\,}(\rho)+t\}g(\rho)~\ge0
&&(g_t(\rho)=0 ~\Leftrightarrow~ 
\rho=0~\textrm{or}~(t=0~\textrm{and}~\rho\le 1/2)), \\
L_t(\rho) &= f_t(\rho)g'_t(\rho)-f'_t(\rho)g_t(\rho)~\ge0
&&(L_t(\rho)=0 ~\Leftrightarrow~
t=0), \\
K_t(\rho) &= f_0(\rho)g'_t(\rho)-f'_0(\rho)g_t(\rho)~\ge0
&&(K_t(\rho)=0 ~\Leftrightarrow~t=0~\textrm{or}~\rho\ge1/2 ), \\
-f'_t(\rho) &=-(1-t)\{b(\rho)f(\rho)\}'-tf'(\rho)~\ge0
&&(f'_t(\rho)=0 ~\Leftrightarrow~f(\rho)=1~\textrm{or}~ (t=0 ~\textrm{and}~ \rho\ge1/2) ),~ \textrm{and} \\
g'_t(\rho) &=(1-t)\{\overline{\,b\,}(\rho)g(\rho)\}'+tg'(\rho))~\ge 0
&&(g_t'(\rho)=0 ~\Leftrightarrow~g(\rho)=1/2~\textrm{or}~ (t=0 ~\textrm{and}~\rho\le1/2))
\end{align*}
Taking further a function $h(\rho)$ such that  
$\mathrm{supp}(h)\subset(0,1)$ and 
$h(1/2)\neq 0$, we put 
\[
\alpha_t= \left\{
\begin{array}{ll}
(1-t)f_0(\rho)\nu+tf_t(\rho)\eta+g_t(\rho)d\theta+(1-t)h(\rho)d\rho
& \textrm{on}\quad N^{2n-1}\times D^2\\
 (1-t^2 )d\theta+t^2\alpha
& \textrm{on}\quad M^{2n+1}\setminus  (N^{2n-1}\times D^2 ).
\end{array} \right.
\]
Note that if the term $(1-t)h(\rho)d\rho$ is dropped, the $1$-form $\alpha_0$ 
at $t=0$ vanishes along the hypersurface $\{\rho=1/2\}$. 
In other words, the term makes the border leaf of the 
dead-end component $\{\rho\le1/2\}$ in the foliation $\F$ defined by $\alpha_0$. 
Note also that the $1$-form $\nu$ on the binding determines the foliation of 
the interior of the dead-end component topologically while it is irrelevant to 
the foliation of the exterior. 
Thus we cut-off $\nu$ by using the bump function $b$ as $bf\nu=f_0\nu$. 
Anyway we see from the following calculations that 
$\alpha_t$ is a contact form for $0<t\le 1$, 
finally $\alpha_1$ coincides with $\alpha$, 
and $\ker\alpha_0$ is actually integrable. 
\begin{itemize}
\item[(i)] On $N^{2n-1}\times D^2$, taking special care of the $2$-form factor 
$d\rho\wedge d\theta$ in the volume form, we have
\begin{align*}
\alpha_t
&=
(1-t)f_0\nu+tf_t\eta+g_td\theta+(1-t)hd\rho,
\\
d\alpha_t
&=
d\rho\wedge\{(1-t)f_0'\nu+tf_t'\eta+g_t'd\theta\}+tf_td\eta,
\\
\alpha_t\wedge (d\alpha_t )^n
&=
\{(1-t)f_0\nu+tf_t\eta\}
\wedge[ng_t'd\rho\wedge d\theta\wedge(tf_td\eta)^{n-1}]
\\
&\qquad +g_td\theta \wedge
[nd\rho\wedge\{(1-t)f_0'\nu+tf_t'\eta\}\wedge(tf_td\eta)^{n-1}]\\
&=
n(t f_t)^{n-1}\{(1-t)K_t\nu+tL_t\eta\}
\wedge(d\eta)^{n-1}\wedge d\rho\wedge d\theta\\
&=\left\{
\begin{array}{ll}
\{(1-t)K_t\nu+tL_t\eta\}\wedge d\rho\wedge d\theta
& (n=1)\\
nt^n(f_t)^{n-1}L_t\eta
\wedge(d\eta)^{n-1}\wedge d\rho\wedge d\theta
& (n>1, \textrm{from Proposition~\ref{prop}}),
\end{array}\right.
\end{align*}
which is positive for $0<t\le 1$, and
\[
\alpha_0\wedge d\alpha_0
=
(f_0\nu+g_0d\theta+ed\rho)\wedge 
d\rho\wedge(f_0'\nu+g_0'd\theta)
=
L_0\nu\wedge d\rho\wedge d\theta
=
0.
\]
\item[(ii)] On $M^{2n+1}\setminus (N^{2n-1}\times D^2)$, we have 
$\alpha_t\wedge (d\alpha_t )^n
=\{ (1-t^2 )d\theta+t^2\alpha\}\wedge  (t^2d\alpha )^n$, which is positive 
for $t>0$ since $d\theta\wedge (d\alpha)^n>0$. On the other hand 
$\alpha_0=d\theta$ obviously defines a foliation.
\end{itemize}
Now we define the reference twisting $\tau$ by perturbing $d\alpha$. 
For a sufficiently small positive constant $\delta$, we put
\[
\tau= \left\{\begin{array}{ll}
d\{f(\rho) \eta+g(\rho) d\theta\}
+\delta\{ h(\rho) d\theta\wedge\nu+\omega\} 
& \textrm{on}\quad N^{2n-1}\times D^2\\
d\alpha +\delta \omega 
& \textrm{on}\quad M^{2n+1}\setminus  (N^{2n-1}\times D^2 ),
\end{array} \right.
\]
where the term $\delta h(\rho) d\theta\wedge\nu$ 
is for making the closed leaf $\{\rho=1/2\}$ almost symplectic. 
Then we see from the following calculation 
that $\alpha_t$ is an $\E\tau$-confoliation form.
\begin{itemize}
\item[(i)] On $N^{2n-1}\times D^2$, 
providing that $0<\delta\le f$ and $\E\in(0,1]$, we have
\begin{align*}
\alpha_t\quad 
&=
(1-t)f_0\nu+tf_t\eta+g_td\theta+(1-t)hd\rho,
\\
d\alpha_t+\E\tau\quad
&=d\rho\wedge
 \{(1-t)f_0'\nu+ (tf_t+\E f )'\eta+ (g_t+\E g )'d\theta
 \}\\
&\qquad +
 (tf_t+\E f )d\eta+\E \delta hd\theta\wedge\nu+\E \delta \omega,
\end{align*}
and
\begin{align*}
\alpha_t\wedge (d\alpha_t+\E\tau )^n
&= \{(1-t)f_0\nu+tf_t\eta \}
\wedge  [n (g_t+\E g )'d\rho\wedge d\theta
\wedge \{ (tf_t+\E f )d\eta
+\E\delta \omega \}^{n-1} ]
\\&\qquad +g_td\theta\wedge
 [nd\rho\wedge \{(1-t)f_0'\nu+ (tf_t+\E f )'\eta \}
\wedge \{ (tf_t+\E f )d\eta+\E\delta \omega \}^{n-1} ]
\\&\qquad +\,(1-t)hd\rho
\wedge [n\E \delta h d\theta\wedge \nu
\wedge  \{ (tf_t+\E f )d\eta+\E\delta \omega \}^{n-1} ]
\\
&=n [(1-t) \{f_0 (g_t+\E g )'-f_0'g_t+\E\delta h^2 \}\nu+
 \{tf_t (g_t+\E g )'- (tf_t+\E f )'g_t \}\eta ]
\\&\qquad\wedge \{ (tf_t+\E f )d\eta+\E\delta \omega \}^{n-1}
\wedge d\rho\wedge d\theta
\\
&=n \{
(1-t) (K_t+\E f_0g'+\E \delta h^2 )\nu
+ (tL_t+\E tf_tg'-\E f'g_t )\eta
 \}
\\&\qquad\wedge
 \{ (tf_t+\E f )d\eta+\E\delta \omega \}^{n-1}
\wedge d\rho\wedge d\theta>0,  
\end{align*}
where the positivity is deduced from the inequality (\ref{ineq}) provided that  
$tL_t+\E (tf_tg'-f'g_t )>0$. Indeed, 
\[
0<\frac{t L_t+\E (t f_t g'-f' g_t )}
{(1-t) (K_t+\E f_0 g' +\E \delta h^2 )+tL_t+\E (tf_tg'-f'g_t )}
\le 1
\]
and
\[
0<\frac{ \{t L_t+\E (t f_t g'-f' g_t ) \}\E \delta}
{ (t f_t+\E f ) \{(1-t)
 (K_t+\E f_0 g' +\E \delta h^2 )+tL_t
+\E (tf_tg'-f'g_t ) \}}
 (\le\frac{\E \delta}{tf_t+\E f}\le
\frac{\delta}{f} )\le 1
\]
can play the roles of $t$ and $\E$ in the inequality (1). 
On the other hand, $tL_t+\E(tf_tg'-f'g_t)=0$ is equivalent to  
$t=0$ and $\rho\le 1/2$. In the case where $n>1$, 
we see from Proposition~\ref{prop} that   
\[
\alpha_0\wedge (d\alpha_0+\E\tau )^n|_{\{\rho\le1/2\}}
=n\E^n (f_0g'+\delta h^2 )\nu\wedge(\delta \omega)^{n-1}
\wedge d\rho\wedge d\theta>0
.\]
It is easy to see that $\alpha_0\wedge (d\alpha_0+\E\tau))|_{\{\rho\le 1/2\}}
=\E\alpha_0\wedge\tau|_{\{\rho\le 1/2\}}>0$ in the case where $n=1$.
\item[(ii)] On the closure of $M^{2n+1}\setminus (N^{2n-1}\times D^2 )$ 
in $M^{2n+1}$, we have
\[\alpha_t\wedge (d\alpha_t+\E\tau )^n
= \{ (1-t^2 )d\theta+t^2\alpha \}
\wedge \{ (t^2+\E )d\alpha+\E\delta\tau \}^{n-1}>0
\]
since $ \{ (1-t^2 )d\theta+t^2\alpha \}
\wedge \{d\alpha+\delta\tau \}^{n-1}>0$ 
holds for sufficiently small $\delta>0$.
\end{itemize}

Next we show that we can take $\omega$ 
which satisfies $d\theta\wedge d\omega=0$ 
in the case where the $2$-form $\tau_N$ is closed 
and extends to a closed $2$-form $\omega_0$ on the page $\Sigma_0$. 
This implies that $\tau|_{T\F}$ is leafwise closed 
while $\tau$ itself can not be closed since it is non-degenerate 
along the closed leaf $\{\rho= 1/2\}$. 
We describe the pages $\Sigma_\theta$ as 
the levels $\Sigma_0\times\{\theta\}$ of a mapping torus 
$\Sigma_0\times[0,2\pi]/(x,2\pi)\sim(\varphi(x),0)$, 
where the restriction of the monodromy map $\varphi$ to the collar 
$\Sigma_0\cap (N^{2n-1}\times D^2 )$ is the identity. 
Take a function $s(\theta): [0,2\pi] \to [0,1]$ which is smoothly 
tangent to $s=0$ and $s=1$ respectively at $\theta=0$ and $\theta=2\pi$. 
Then the $2$-form $\omega=(1-s(\theta))\omega_0+s(\theta)\varphi^*\omega_0$ 
satisfies $d\theta\wedge d\omega=0$. Further we may assume that 
the restriction of $\omega_0$ to the collar is the pull-back of $\tau_N$. 
Then the restriction of $\omega$ to the tube $N^{2n-1}\times D^2$ is 
also the pull-back of $\tau_N$. This completes the proof. \qed  

\begin{remark}
The leaf space $M^{2n+1}/\F$ contains $S^1$ presenting 
the closed transversal $\{*\}\times S^1 \subset N^{2n-1}\times D^2$. 
It also contains the points representing the closed leaves $\{\rho=1/2\}$. 
The complement of their union can be identified with the leaf space of the 
Riemannian foliation on $N^{2n-1}$ defined by the closed $1$-form $\nu$. 
If $\{\rho=0\}$ is connected and $\nu$ defines a fiber bundle, 
then the leaf space of $\F$ is homeomorphic to that of the Reeb foliation of $S^3$. 
In contrast to the result of Meigniez \cite{Meigniez} on the flexibility of 
codimension one smooth foliations of closed manifolds of dimension $>3$, 
Mart\'inez Torres \cite{DMT} proved that if a leafwise symplectic structure 
of a codimension one foliation is the restriction of a closed $2$-form on the manifold, 
the leaf space is homeomorphic to that of a taut foliation of a $3$-manifold. 
Thus it is natural to seek a wider class of leafwise symplectic foliations 
which inherits another three dimensional property such as the Novikov closed leaf theorem. 
\end{remark}

\section{Application}
As an application of Theorem~\ref{thm}, we construct 
a family of leafwise symplectic foliations of $S^4\times S^1$. 
\subsection{The Furukawa-Kasuya construction}
Let $z_i=r_i e^{\sqrt{-1}\varphi_i}$ be the $i$-th 
coordinate of $\C^3$. 
The $5$-sphere $S^5= \{r_1^2+r_2^2+r_3^2=1 \}\subset \C^3$ 
carries the standard contact form 
$\alpha=r_1^2d\varphi_1+r_2^2d\varphi_2+r_3^2d\varphi_3$,  
the standard open-book decomposition  
$\varphi_1:S^5\setminus\{r_1=0\}\to S^1$ which supports $\ker\alpha$, 
and the singular torus fibration 
$\textrm{pr}=(r_1, r_2, r_3): S^5\to \Delta\subset \R^3$ to the spherical $2$-simplex  
$\Delta$ which is the intersection of the unit sphere and the positive orthant. 
For each positive integer $p$, we will customize the contact form 
\[
\alpha_p=\frac{1}
{2r_1^2+u_p (r_1^2 ) (1-r_1^2 )}~\alpha
\]
and the map
\[
w_p=v_p (|z_1|^2)z_1^p +z_2 z_3=v_p (r_1^2 ) r_1^p e^{p\sqrt{-1}\varphi_1}
+r_2r_3 e^{\sqrt{-1}(\varphi_2+\varphi_3)}: S^5 \to \C
\] 
by using non-negative functions $u_p(x)$ and $v_p(x)$ of $x\in [0,1]$ 
which will be defined in the sequel. Furukawa and Kasuya \cite{Kas} 
following the author's construction 
of a contact submanifold of $S^5$ in \cite{MoriCR} arranged these functions 
so that the zero set $N_p^3=\{w_p=0\}\subset S^5$ of $w_p$ inherits 
a contact form $\eta_p$ and a torus fibration $\textrm{pr}_p: N_p^3\to a\subset \Delta$ 
over an arc $a$ connecting $(0,1,0)$ to $(0,0,1)$ through $\textrm{int}~\Delta$ 
such that $\ker\eta_p$ is a propeller contact structure with respect to $\textrm{pr}_p$, 
and therefore isomorphic to the standard contact structure on the Lens space $L(p,p-1)$. 
Our task is to improve their construction so that further  
$\theta_p=\arg w_p:S^5\setminus N_p^3\to S^1$ 
gives a supporting open-book decomposition which is different from $\varphi_1$. 
Namely we show the additional inequality \begin{equation}
d\theta_p\wedge (d\alpha_p)^2>0
\end{equation}
under a careful choice of the functions $u_p$ and $v_p$.
 
First we assume that $v_p(x)=0$ for $0\le x\le 1$. 
Put $\rho_i=r_i^2$ ($i=1,2,3$) and write the inequality $(2)$ as  
\begin{align*}
d(\rho_1+\rho_2+\rho_3)\wedge d\theta_p\wedge (d\alpha_p)^2
&=d(\rho_1+\rho_2+\rho_3)\wedge 
d(\varphi_2+\varphi_3)\wedge (d\alpha_p)^2\\
&=\frac{4 \{u_p(\rho_1)-\rho_1(1-\rho_1)u_p'(\rho_1) \}}
{\{2\rho_1+u_p(\rho_1)(1-\rho_1)\}^3}
d\rho_1\wedge d\varphi_1
\wedge d\rho_2\wedge d\varphi_2
\wedge d\rho_3\wedge d\varphi_3>0
\end{align*}
in $\C^3=\R^6$. 
We look at the foliation 
$\displaystyle \left\{y=cx/(1-x)\right\}_{c>0}$ which presents the general solution 
of the differential equation $y-x(1-x)y'=0$. 
If the graph $\displaystyle \{(x,y) \mid y=u_p(x), 0\le x<1 \}$ 
is transverse to this foliation, 
the sign of the above $6$-form does not change. 
Hence we choose the function $u_p(x)$ so that 
\begin{itemize}
\item[(i)] $u_p(x)=2$ holds for $\displaystyle 0 \le x\le \frac{1}{2}$,
\item[(ii)] $u_p(x)=p$ holds for $\displaystyle \frac{p+1}{p+2}<x<1$, and
\item[(iii)] the graph $ \left\{(x,y) \mid y=u_p(x), 0\le x<1 \right\}$ is transverse to 
the foliation $\left\{y=cx/(1-x)\right\}_{c>0}$ 
\end{itemize}
This is possible even when $p>2$ since 
the leaf $\displaystyle\left\{y=2x/(1-x)\right\}$ through 
the point $(1/2,2)$ reaches the point $\displaystyle\left(p/(p+2\right),p)$, and $p/(p+2)<(p+1)/(p+2)$ holds. 
This implies the inequality $(2)$ except along 
the singular binding $\{r_2r_3=0\}$ of the singular open-book decomposition 
$\varphi_2+\varphi_3: S^5\setminus  \{r_2r_3=0 \}\to S^1$. Here the singularity is 
the circle $\{r_1=1\}$ along which the smooth binding meets itself transversely. 

Next we perturb the function $v_p(x)\equiv 0$ on the interval 
$\displaystyle \left((p+1)/(p+2),~1\right]$ so that it is 
smoothly tangent to a positive constant at $x=1$. 
Then we see that the vector field $X_p=2\partial_{\varphi_1}
+p\partial_{\varphi_2}+p\partial_{\varphi_3}$  
is the Reeb vector field of $\alpha_p$ on the neighborhood 
$\displaystyle\left \{(p+1)/(p+2)<r_1^2\le 1 \right\}$ of the singular binding $\{r_2r_3=0\}$, and it satisfies 
$d(\Re (w_p))(X_p)=-2p\Im (w_p)$ and $d(\Im (w_p))(X_p)=2p\Re (w_p)$. 
This implies that $d\theta_p (X_p )=2p>0$ 
hence the inequality $(2)$ holds for the regular open-book decomposition $\{\theta_p=\textrm{const}\}$.

In the case where $p>1$, if we replace the function 
$v_p (r_1^2 )$ with 
a positive constant, the open-book decomposition becomes 
the Milnor fibration of the isolated singularity $(0,0,0)$ 
of the complex surface $z_1^p=z_2z_3$. Further if $p>2$, 
it approximate the singular open-book decomposition $\varphi_2+\varphi_3$ 
of a small hypersphere. 

\subsection{Leafwise symplectic foliations of $S^4\times S^1$}
We take a copy of $S^5$ together with all of the above objects
and use notations with tildes to denote the copy. 
Then 
$\widetilde{S}^5= \{\widetilde{r}_1^2+\widetilde{r}_2^2+\widetilde{r}_3^2=1 \}$ 
carries the open-book decomposition 
$\widetilde{\theta}_{\widetilde{p}}:
\widetilde{S}^5\setminus \widetilde{N}^3_{\widetilde p}\to S^1$
for any positive integer $\widetilde{p}$ which may differ from $p$. 
We remove the binding $ \{r_1=0, r_2^2+r_3^2=1  \}$ from
the trivial open-book decompositions of $S^5$, also remove 
the binding $ \{\widetilde{r}_1=0, \widetilde{r}_2^2+\widetilde{r}_3^2=1 \}$
from that of $\widetilde{S}^5$, and 
glue the rests along $S^3\times S^1$ 
together by the identification $r_2\sim\widetilde{r}_2$, $r_3\sim\widetilde{r}_3$, 
$\varphi_2\sim\widetilde{\varphi}_2$, $\varphi_3\sim\widetilde{\varphi}_3$ 
and $\varphi_1\sim-\widetilde{\varphi}_1$. 
Precisely, we extend the angular coordinate $\varphi_1$ (resp. $\widetilde{\varphi}_1)$ 
to the boundary of the exterior of the 
unknot exterior $\{r_1=0\}$ (resp. 
$\{\widetilde{r}_1=0\}$), and replace the normal coordinate $r_1$ 
(resp. $\widetilde r_1$) of the boundary with $h=r_1^2$ (resp. 
with $-h=\widetilde{r}_1^2$). 
Then the contact forms $\alpha_p$ and $\widetilde{\alpha}_{\widetilde{p}}$ 
match up to define a contact form $\alpha_{p,\widetilde{p}}$ on $S^4\times S^1$. 
Note that one can also obtain 
this contact structure just by smoothing the boundary of 
the symplectic product $A\times B^4$ of the annulus $A$ and the $4$-ball. 
Indeed one can break up the product $A\times B^4$ of the annulus $A$ and 
the $4$-ball into 
two copies of (hemisphere)$\times D^2\approx D^2\times B^4$ 
by shrinking the core of the first factor $A$.  
The following are the advantages of our precise construction. 
\begin{itemize}
\item[(i)] The page fibrations 
$ \{\theta_p=\mathrm{const} \}$ and 
$ \{\widetilde{\theta}_{\widetilde{p}}=\mathrm{const} \}$ 
match up to define a \textit{smooth} open-book 
decomposition $\theta_{p,\widetilde{p}}: S^4\times S^1
\setminus N^3_{p,\widetilde{p}}\to S^1$ 
which supports the contact structure $\ker\alpha_{p,\widetilde{p}}$.  
We put $M^5=S^4\times S^1$, 
$\theta=\theta_{p,\widetilde{p}}$, and $\alpha=\alpha_{p,\widetilde{p}}$ 
in Theorem~\ref{thm}. 
\item[(ii)] We can take the closed $2$-form $\omega$ in Theorem~\ref{thm}
by restricting 
$d\varphi_1\wedge (d\varphi_2-d\varphi_3 )$ and 
$d\widetilde{\varphi}_1\wedge
 (d\widetilde{\varphi}_2-d\widetilde{\varphi}_3 )$ 
to the page $\Sigma_0= \{\theta_{p,\widetilde{p}}=0 \}$ and extending it 
to the boundary $\partial\Sigma_0=N_{p,\widetilde{p}}^3$. 
Note that if we fix the ratio $r= [v_p (r_1^2 )r_1^p:r_2r_3 ]$, 
the function $\theta_p$ depends only on the 
two variables $\varphi_1$ and $\varphi_2+\varphi_3$. 
It becomes independent of $\varphi_2+\varphi_3$ as $r\to[1:0]$. 
Thus $\omega$ vanishes along $ \{r_2r_3=0 \}\cap \Sigma_0$. 
\item[(iii)] We can take the closed $1$-form $\nu$ 
in Theorem~\ref{thm} as the pull-back of a positive $1$-form on the circle 
$a\cup \widetilde{a}$ on the spherical lune $\Delta\cup \widetilde{\Delta}$ 
under the $T^2$-bundle projection. 
Let $T_r$ denote the fiber at $r\in a\cup\widetilde{a}$.
\end{itemize}

Hereafter we omit the description in $\{h\le 0\}$ for simplicity. 
In the case where $r\in a\setminus \partial\Delta$, 
the restriction $\omega|_{T_r}$ is obviously non-degenerate 
since $T_r$ is the section $\{p\varphi_1=\varphi_2+\varphi_3+\pi\}$ 
of $T^3\ni(\varphi_1,\varphi_2,\varphi_3)$. 
In the other case, we have 
$\theta_p=\varphi_2+\varphi_3=0$ near $T_r$. If 
$r$ is on the edge $\{y=0\}\cap \Delta$, then
we can write $\omega=2d\varphi_1\wedge d\varphi_3$, 
else $r\in\{z=0\}\cap\Delta$ and $\omega=2d\varphi_1\wedge d\varphi_2$. 
Thus $\omega|_{T_r}$ is also non-degenerate. 
Therefore Theorem~\ref{thm} provides 
a leafwise symplectic foliation $\F_{p,\widetilde{p}}$ on $S^4\times S^1$
which can be deformed into the contact structure $\ker\alpha_{p,\widetilde{p}}$.

\begin{remark}
Since $N_p$ is $L(p,p-1)$, we can see that the binding 
$N_{p,\widetilde{p}}\subset  (S^4\times S^1,\ker\alpha_{p,\widetilde{p}} )$ 
is the contact $Nil$-manifold arising as the mapping torus 
of the action on $T^2$ by the parabolic element 
$\displaystyle \left[
\begin{array}{cc}
1&0\\p+\widetilde{p} &1\end{array} \right]$ of $SL(2,\Z)$. 
You can see it directly tracing the isotopy of the 
vector fields $(1,0,0)$ and $(0,1,0)$ on $T^3$ 
starting at the middle of the left arc of the spherical lune 
$\Delta\cup\widetilde{\Delta}$ and moving 
along $a\cup \widetilde{a}$ once clockwise with respect to the height $h$. 
They respectively jump up to $(1,0,p)$ and $(0,1,-1)$, and fall down to 
$(1,0,p)$ and $(0,0,-1)$ since temporarily $p\varphi_1=\varphi_2+\varphi_3+\pi$ 
holds on the way along $a$. 
Then they respectively jump up to $(1,-p-\widetilde{p},p)$ and $(0,1,-1)$, 
and fall down to $(1,-p-\widetilde{p},p)$ and $(0,1,0)$ since temporarily 
$-\widetilde{p}\varphi_1=\varphi_2+\varphi_3+\pi$ 
holds on the way along $\widetilde{a}$. 
Thus the monodromy $\mu\in SL(2,\Z)$ satisfies 
$\displaystyle 
\mu \left[ \left[\begin{array}{c}1\\-p-\widetilde{p}\end{array} \right]
 \left[\begin{array}{c}0\\1\end{array} \right] \right]=I.
$ 
Although $\F_{p,\widetilde{p}}$ is clearly different from $\F_{q,\widetilde{q}}$ 
provided $p+\widetilde{p}\neq q+\widetilde{q}$, 
the author do not know whether the sum $p+\widetilde{p}$ 
determines $\F_{p,\widetilde{p}}$ or not.
\end{remark}

\subsection{Towards complex geometric understanding}
It is natural to seek a complex geometric model of 
$N_{p,\widetilde{p}}$ in the case where both $p$ and $\widetilde{p}$ 
are greater than $1$. 
In this task it would be suggestive that 
among simple elliptic singularities only the three  
contact $Nil$-manifolds associated with
$\displaystyle \left[
\begin{array}{cc}
1&0\\1&1\end{array} \right]$, 
$\displaystyle \left[
\begin{array}{cc}
1&0\\2&1\end{array} \right]$, and 
$\displaystyle \left[
\begin{array}{cc}
1&0\\3&1\end{array} \right]$ 
can be realized as the link of hypersurface singularities 
(see for example \cite{Saito}). 
Indeed for any positive integers 
$p,q,r$ with $1/p+1/q+1/r\le 1$, 
we have the Furukawa-Kasuya supporting open-book decomposition of $(S^5,\ker\alpha)$ 
by pages 
\[\Sigma_\theta= \{
\arg (v_{p-1} (r_1^2 ) r_1^p e^{p\sqrt{-1}\varphi_1}
+v_{q-1} (r_2^2 ) r_2^q e^{q\sqrt{-1}\varphi_2}
+v_{r-1} (r_3^2 ) r_3^r e^{r\sqrt{-1}\varphi_3}
+r_1r_2r_3 e^{\sqrt{-1}(\varphi_1+\varphi_2+\varphi_3)} )=\theta \}, 
\] 
which approximates the Milnor fibration 
of the singularity $(0,0,0)$ of the surface $z_1^p+z_2^q+z_3^r+z_1z_2z_3=0$ 
(see Kasuya \cite{Kas}). 
Here the binding $\partial\Sigma_0$ is a
contact $Nil$/$Sol$-manifold arising as the mapping 
torus of the parabolic/hyperbolic element 
$\displaystyle \left[
\begin{array}{cc}
r-1&-1\\1&0\end{array} \right] \left[
\begin{array}{cc}
q-1&-1\\1&0\end{array} \right] \left[
\begin{array}{cc}
p-1&-1\\1&0\end{array} \right]\in SL(2,\Z)$. 
If it is parabolic, it is conjugate to $\displaystyle \left[\begin{array}{cc}
1&0\\k&1\end{array} \right]$ for $k=1$, $2$, or $3$ 
according to $\{p,q,r\}=\{2,3,6\}$, $\{2,4,4\}$, or $\{3,3,3\}$ (as a multiset). 
We can also take the closed $2$-form $\omega$ 
in Theorem~\ref{thm} so that $\omega$ is 
expressed in the domains $\{r_i^2>1/2\}$ 
as $d\varphi_i\wedge (d\varphi_{i+1}-d\varphi_{i+2})$ ($i\in\Z_3
$), 
and elsewhere as $2d\varphi_2\wedge d\varphi_3
=2d\varphi_3\wedge d\varphi_1=2d\varphi_1\wedge d\varphi_2$. 
This reproduces Mitsumatsu's leafwise symplectic foliations, 
as well as relates them to the standard contact structure.

\begin{flushright}
Department of Mathematics, Osaka Dental University\\ 
8-1 Kuzuha-hanazono-cho, Hirakata, Osaka 573-1121 JAPAN\\
mori-a@cc.osaka-dent.ac.jp, ka-mori@ares.eonet.ne.jp
\end{flushright}
\end{document}